\documentclass[12pt]{article}
\usepackage[centertags]{amsmath}
\usepackage{amsfonts}
\usepackage{amssymb}
\usepackage{amsthm}
\usepackage{newlfont}
\usepackage{graphicx}
\newtheorem{Theorem}{Theorem}[section]

\newtheorem{Lemma}{Lemma}[section]
\newtheorem{Remark}{Remark}[section]

\theoremstyle{definition}

\theoremstyle{remark}

\numberwithin{equation}{section}

\begin{document}

\title{On the asymptotic stability of a rational multi-parameter first order difference equation}
\author{M. Shojaei\footnote{ Corresponding author, \ \ \textit{E-mail
address}: m\_shojaeiarani@aut.ac.ir,}}
\date{}
\maketitle
\begin{center}
\textit{\scriptsize Department of Applied
Mathematics, Faculty of Mathematics and Computer Science,
Amirkabir University of Technology, No. 424, Hafez Ave, Tehran,
Iran.}\\[0pt]
\end{center}
\begin{abstract}
In this part we study the dynamics of the following rational
multi-parameter first order difference equation
$$x_{n+1}=\frac{ax_{n}^3+bx_{n}^2+cx_{n}+d}{x_{n}^3},\ \ \ x_{0}\in \Bbb{R}^+$$
where the parameters $a,b,d$ together with the initial condition
$x_{0}$ are positive while the parameter $c$ could accept some
negative values. We investigate the equilibria and 2-cycles of this
equation and analyze qualitative and asymptotic behavior of
it's solutions such as convergence to an equilibrium or to a
2-cycle.
\end{abstract}
\section{Introduction}
Difference equations and discrete dynamical systems appear as both
the discrete analogs of differential and delay differential
equations in which time and space are not continuous and as direct
mathematical models of diverse phenomena such as biology (see
\cite{BC,B}), medical sciences (see \cite{MK}), economics (see
\cite{IS,LM}), military sciences (see \cite{Epstein,Sedaghat1}), and
so forth.

One of the most practical classes of nonlinear difference equations
(for nonlinear difference equations see
\cite{AXS,Sedaghat2,SHS,SSA})) are rational difference equations
that are simply the ratio of two polynomials. Most of the work about
rational type difference equations treat the case when the numerator
and denominator are both linear polynomials (see the monograph of
Kulenovic and Ladas \cite{KL}). Also there are some works about
rational difference equations of order two with quadratic terms in
numerator or denominator (see \cite{CLNR,CK,CKS,DKMOS,DKMOS1,GJKL}).

In this paper we analyze the global behavior of the following
rational difference equations
\begin{equation}\label{formula1}
x_{n+1}=\frac{ax_{n}^3+bx_{n}^2+cx_{n}+d}{x_{n}^3},
\end{equation}\\
where the parameters $a,b,d$ together with the initial value $x_{0}$
are positive while the parameter $c$ could accept some negative
values. We study the equilibria and 2-cycles of Eq.(\ref{formula1})
and show that in most cases every positive solution of this equation
converges to either an equilibrium or a 2-cycle.

Suppose that the sequence $\{x_{n}\}$ is a solution of the first
order difference equation  $x_{n+1}=f(x_{n})$, where $f:I\rightarrow
I$, $I\subseteq \Bbb{R}$. The solution $\overline{x}$ of the
equation $x=f(x)$ is called an equilibrium (or fixed) point of this
equation. The equilibrium point $\overline{x}$ is said to be stable
if for every $\epsilon
>0$ there exists $\delta >0$ such that for any initial condition
$x_{0}\in I$ with $|x_{0}-\overline{x}|<\delta $, the iterates
$x_{n}$ satisfy $|x_{n}-\overline{x}|<\epsilon $ for all $n\in \Bbb{N}$.
$\overline{x}$ is said to be attracting if for all $x_{0}\in I,\lim
_{n\rightarrow \infty }x_{n}=\overline{x}$. $\overline{x}$ is
asymptotically stable if it is both stable and attractor. The point
$p\in I$ is called a period two solution if $f^{2}(p)=p$, where
$f^{2}$ denotes the second iterate of $f$. If moreover, $f(p)\neq p$
then $p$ is called a prime period two solution. In this case the
point ($p$,$f(p)$) is called a $2$-cycle of $f$. An interval
$I\subseteq \Bbb{R}$ is called invariant under $f$ if $f(I)\subseteq
I$. For higher order difference equations similar notions could be
defined (see \cite{Elaydi}).
\section{The parameter $c$ and equilibria}
\begin{Theorem}\label{t1} Let $\phi (x)=\frac{ax^3+bx^2+cx+d}{x^3},x>0$ which is the right hand side of Eq.(\ref{formula1}).
\begin{description}
    \item[\it{(a)}] The following cubic polynomial has a unique negative zero
    $$Q(x)=(4a)x^3-(b^2)x^2-(18abd)x+27a^2d^2+4db^3.$$
    \item[\it{(b)}] Assume that $c_{-}$ is the unique negative zero of the cubic polynomial in $(a)$
    . If $c>c_{-}$ then nonpositive iterations of Eq.(\ref{formula1})
    never occur.
    \item[\it{(c)}] Let $c^{*}=-\sqrt{3bd}$. Then $c_{-}<c^{*}$.
    \item[\it{(d)}] If $c\geq c^{*}$ then $\phi$ is decreasing and
    has a unique equilibrium. If $c_{-}<c<c^{*}$ then $\phi$ has a minimum  point $x_{m}$ and
    a maximum point $x_{M}$ where
    $$0<x_{m}=\frac{-c-\sqrt{c^2-3bd}}{b}<
    x_{M}=\frac{-c+\sqrt{c^2-3bd}}{b}.$$
\end{description}
\end{Theorem}
Proof. ($a$) The facts that $Q(0)>0,Q'(0)<0,Q(x)\rightarrow
-\infty $ as $x\rightarrow -\infty$ together with the fact that $Q$
is a cubic polynomial prove ($a$).

($b$) Define $F(x)=ax^3+bx^2+cx+d,x> 0$ which simply is the
numerator of $\phi$. If $F(x)>0$ then Eq.(\ref{formula1}) could not
accept nonpositive iterations. Note that $F(x)>0$ if $c\geq 0$. Now,
we want to determine the negative values of $c$ such that $F(x)>0$.
Thus, assume that $c<0$. Since $F$ decreases as the parameter $c$
decreases then it's evident that for
a special negative value of $c$ there will exist a point $x^{*}$ such that\\
$$F(x^{*})=F'(x^{*})=0.$$

Some calculations show that the equation $F^{'}(x^{*})=0$ implies
that $bx^{*^{2}}=-\frac{2b^{2}}{3a}x^{*}-\frac{bc}{3a}$ and
$ax^{*^{3}}=\frac{4b^{2}}{9a}x^{*}-\frac{2bc}{3a}$. Replace
$bx^{*^{2}}$ and $ax^{*^{3}}$ into the equation $F(x^{*})=0$ to
obtain
$$x^{*}=\frac{bc-9ad}{6ac-2b^{2}}.$$

Note that $x^{*}>0$ since $c<0$. Now, in order to obtain the value
of the parameter $c$ which corresponds to $x^{*}$ replace $x^{*}$ in
the equation $F'(x^{*})=0$ and solve for $c$ to obtain
$$(4a)c^3-(b^2)c^2-(18abd)c+27a^2d^2+4db^3=0,$$
which simply is the equation $Q(c)=0$ and by $(a)$ this equation has
a unique negative zero, namely $c_{-}$. Therefore, in order to avoid
nonpositive iterations for Eq.(\ref{formula1}) the parameter $c$
should be greater than $c_{-}$.

($c$) Note that $Q(c^{*})=6abd\sqrt{3bd}+db^3+27a^2d^2>0$ ,
$c^{*}<0$, and $c_{-}$ is the unique negative zero of $Q$. These
facts prove $(c)$.

($d$) The equation $\phi '(x)=0$ is equivalent to the equation
$bx^2+2cx+3d=0$ with the determinant $\Delta =4(c^2-3bd)$. Thus if
$c^{*}\leq c  \leq -c^{*}$, then $\phi $ is decreasing on $(0,\infty
)$. Also for $c>-c^{*}$ $\phi $ has no positive extremum since the
summation and the product of the extremum points are negative and
positive respectively. Therefore, if $c\geq -c^{*}$ then again $\phi
$ is decreasing on $(0,\infty )$. So, $\phi $ is decreasing on
$(0,\infty )$ if $c\geq c^{*}$. Hence, $\phi $ has a unique
equilibrium when $c\geq c^{*}$. Now, assume that $c_{-}<c<c^{*}$. In
this case $\phi $ has a minimum point $x_{m}=(-c-\sqrt{c^2-3bd})/b$
and a maximum point $x_{M}=(-c+\sqrt{c^2-3bd})/b$, both positive. The proof is complete.\\

 In order to avoid nonpositive iterations we assume that $c>c_{-}$, hereafter.
The next theorem deals with the number of equilibria when $c\in
(c_{-},c^{*})$.\\

\begin{Theorem}\label{t2} Suppose that $c\in
(c_{-},c^{*})$ and let $P(t)=t^4-at^3-bt^2-ct-d, \ t>0$. Suppose
also that $c_{b}$ is the unique negative root of the following
quadratic polynomial
$$H(x)=108x^2+(108ab+27a^3)x-9a^2b^2-32b^3.$$
\begin{description}
\item[\it{(a)}] If $c=c_{b}$ then $P'$ touches the horizontal axis at
$t^{*}=-\frac{6c_{b}+ab}{3a^2+8b}>0$.
\item[\it{(b)}] If $c_{b}\geq \frac{b^2-12d}{3a}$ then $\phi $ has
 a unique positive equilibrium for all $c\in (c_{-},c^{*})$.
\item[\it{(c)}] If $c_{b}< \frac{b^2-12d}{3a}$ then by increasing
the parameter $c$ from $c_{b}$ $P$ touches the horizontal axis at a
minimum point $t_{m}$ (at first) and a maximum point $t_{M}$ (next)
with $0<t_{M}<t_{m}$. Assume that $c_{m}$ and $c_{M}$ are the values
of parameter $c$ which correspond to $t_{m}$ and $t_{M}$
respectively. Then
$$c_{b}<c_{m}<c_{M}<c^{*}, \ \ \ c_{-}<c_{M}.$$
In this case there are two subcases as follow:
\begin{description}
    \item[\it{($c_{1}$)}] $c_{-}\leq c_{m}$: If $c\in (c_{m},c_{M})$ then $\phi $ has
    three positive equilibria. For $c=c_{m}$, $\phi $ has two positive equilibria. One of them is $t_{m}$ with
    the fact that $\phi $ is tangent to the 45 degree line at $t_{m}$ and $t_{m}$ is greater than the other equilibrium.
    For $c=c_{M}$, $\phi $ has again two positive equilibria. One of them is $t_{M}$ with
    the fact that $\phi $ is tangent to the 45 degree line at $t_{M}$ and $t_{M}$ is lower than the other
    equilibrium. Finally, if $c\in (c_{-},c_{m})\cup (c_{M},c^{*})$ then $\phi $
    has a unique positive equilibrium.
    \item[\it{($c_{2}$)}] $c_{m}<c_{-}$;\ If $c\in
    (c_{-},c_{M})$ then $\phi $ has three positive equilibria. For $c=c_{M}$
    $\phi $ has two positive equilibria. One of them is $t_{M}$ with
    the fact that $\phi $ is tangent to the 45 degree line at $t_{M}$ and $t_{M}$ is lower than the
    other equilibrium. Finally, if $c\in (c_{M},c^{*})$ then $\phi $ has a unique
    positive equilibrium.
\end{description}
\end{description}
\end{Theorem}
Proof. (a) With an analysis precisely similar to what was applied in
Theorem 1(b) we obtain that $P'$ touches the horizontal axis at
$t^{*}$ when $c=c_{b}$. It remains to show that $t^{*}>0$. A little
algebra shows that
\begin{equation*}
 H(-ab/6)=-24a^2b^2-9/2a^4b-32b^3<0,
\end{equation*}
this together with the fact that $H(c_{b})=0$ imply that
$c_{b}<-ab/6$, hence $t^{*}>0$.

(b) At first note that the roots of $p$ are simply the fixed points
of $\phi $. Note also that $P'' (t^{*})=P'(t^{*})=0,\ \partial P'
(t)/\partial c =-1<0$ and $P'$ is decreasing on $(0,t^{*})$ and
increasing on $(t^{*},\infty )$ for any $c$. Also, $P' (0)=-c>0$ for
all $c<0$. These facts imply that $P'$ has no positive root for all
$c< c_{b}$ and two positive roots for all $c\in (c_{b},0)$. Hence
$P$ has no positive extremum for all $c\leq c_{b}$ and two positive
extremums for all $c\in (c_{b},0)$. Therefore, if $c\leq c_{b}$ then
$P$ is increasing on $(0,\infty )$ and will have a unique positive
root (note that by the intermediate value theorem $P$ has at least
one positive root). Now, assume that $c=c_{b}$. By the relation
$c_{b}^2=1/108(9a^2b^2+32b^3-(108ab+27a^3)c_{b})$ and by some
algebra one can write
\begin{equation*}
    P(t^{*})=-1/4ac_{b}+1/12b^2-d.
\end{equation*}

so $P(t^{*})\leq 0$ since $c_{b}\geq (b^2-12d)/3a$. This fact
together with the fact that $\partial P(t)/\partial c=-t<0$ imply
that $P$ has a unique positive root when $c>c_{b}$, too. Thus $P$
has a unique positive root for all $c<0$, hence for $c\in
(c_{-},c^{*})$. So $\phi $ has a unique positive equilibrium for
$c\in (c_{-},c^{*})$.

(c) Since $P'$ is decreasing on $(0,t^{*})$ and increasing on
$(t^{*},\infty )$ then the local maximum of $p$ is lower than it's
local minimum for all $c\in (c_{b},0)$. Therefore, according to the
facts that $\partial P(t)/\partial c =-t<0$ and for $c=c_{b}$
$P(t^{*})>0$ (since $c_{b}< (b^2-12d)/3a$) it is evident that as the
parameter $c$ increases from $c_{b}$, $P$ will touch the horizontal
axis two times, at first in a local minimum point (namely $t_{m}$
when $c=c_{m}$) and then in a local maximum point (namely $t_{M}$
when $c=c_{M}$). Thus $c_{b} <c_{m}<c_{M}$.

 Now, suppose that $c<c_{m}$. Then $P$ has a unique root and by
 the intermediate value theorem this root is less than the maximum
 point of $P$. After that, suppose that $c=c_{m}$. Then $P$ touches
 the horizontal axis at $t_{m}$. It is easy to show that this is equivalent
 to this fact that $\phi $ is tangent to the 45 degree line at $t_{m}$. In
 this case $P$ has another root which by the intermediate value
 theorem is less than the maximum point of $P$ and therefore is less
 than $t_{m}$. Next, assume that $c\in (c_{m},c_{M})$. Then, using
 intermediate value theorem, $P$ has a unique root less that the maximum
 point of $P$, a unique root between the maximum and minimum points
 of $P$, and finally a unique root greater than the minimum point of
 $P$. So $P$ has three roots in this case. Similar to the case
 $c=c_{m}$,  in the case $c=c_{M}$, $P$ has two roots with $t_{M}$
 as the lower root. Finally, similar to the case $c<c_{m}$, in the case $c>c_{M}$ $P$ has
 a unique root greater than the minimum point of $P.$

It's easy to show that $\phi (x_{m})\geq x_{m}>0$ for $c\geq c_{M}$,
where $x_{m}$ is the local minimum of $\phi $. Also, by Theorem
$1(b)$ $\phi (x_{m})(c-c_{-})>0 $ for $c\neq c_{-}$ and $\phi
(x_{m})=0$ for $c=c_{-}$. These facts together with the fact that
$\partial \phi /\partial c>0$ simply prove that $c_{-}<c_{M}$. The
other details of the proof of $(c_{1})$ and $(c_{2})$ are
straightforward
and will be omitted. The proof is complete.\\

\begin{Theorem}\label{t3}
\begin{description}
    \item[\it{(a)}] If $\phi $ has a unique positive equilibrium $\overline{t}$
    then
    $$(\phi (t)-t)(t-\overline{t})<0, \ \ \ t>0,\ t\neq \overline{t}. $$
    \item[\it{(b)}] If $\phi $ has two equilibria
    $\overline{t}_{1}<\overline{t}_{2}$ then for $c=c_{M}$
         $$(\phi (t)-t)(t-\overline{t}_{2})<0, \ \ \ t>0,\ t\neq \overline{t}_{1},\overline{t}_{2},$$
    and for $c=c_{m}$
         $$(\phi (t)-t)(t-\overline{t}_{1})<0, \ \ \ t>0,\ t\neq \overline{t}_{1},\overline{t}_{2},$$
    where $c_{m}$ and $c_{M}$ have introduced in Theorem 2.
    \item[\it{(c)}] If $\phi $ has three equilibria
    $\overline{t}_{1}<\overline{t}_{2}<\overline{t}_{3}$ then
    $$(\phi (t)-t)(t-\overline{t}_{1})(t-\overline{t}_{2})(t-\overline{t}_{3})<0,\ \ \ t>0,\ t\neq
    \overline{t}_{1},\overline{t}_{2},\overline{t}_{3}.$$
\end{description}
\end{Theorem}
Proof. The proof is clear (using Theorem 2) and will be omitted.\\
\section{Two cycles and convergence}
\begin{Lemma}\label{Lemma1} Assume that $p$ is a prime period two solution of $\phi $.
 Then $p$ is a zero of the following sixth order polynomial\\
 $$G(x)=a^3x^6+(2a^2b-ac-d)x^5+(ab^2-ad-bc+2a^2c)x^4+(2a^2d+2abc-c^2-bd)x^3$$
 $$+(ac^2+2abd-2cd)x^2+(2acd-d^2)x+ad^2,$$
 Therefore $\phi $ has at most three $2$-cycles.
\end{Lemma}
Proof. Some algebra shows that\\
\begin{equation}\label{f3}
\phi ^{2}(t)-t=\frac{t^3}{(at^3+bt^2+ct+d)^3}\ (\phi (t)-t)G(t).
\end{equation}

Since $p$ is a prime period two solution then $\phi (p)\neq p,\ \phi
^{2}(p)=p$. Hence Eq.(\ref{f3}) implies that $G(p)=0$. Now note that
if $(p,q)$ be a $2$-cycle then both $p$ and $q$ are roots of $G$.
This simply shows that the number of $2$-cycles of $\phi $ is at
most three. The proof is complete.

\begin{Lemma}\label{Lemma2} Assume that $c\in(c_{-},c^{*})$ and
$\overline{t}$ is an equilibrium of $\phi $ with $x_{M}\leq
\overline{t}$. Define
$$G_{1}=-dx_{M}^{4}+(ad-bc)x_{M}^{3}-(c^{2}+bd)x_{M}^{2}-2cdx_{M}-d^{2},$$
and
$$G_{2}=a^{2}x_{M}^{6}+(2ab-c)x_{M}^{5}+(b^{2}-2d+2ac)x_{M}^{4}+(2ad+2bc)x_{M}^{3}+(c^{2}+2bd)x_{M}^{2}+2cdx_{M}+d^{2}.$$
\begin{description}
    \item[\it{(a)}] Both $G_{1}$ and $G_{2}$ are positive.
    \item[\it{(b)}] $\phi ^{2}(x_{M})>x_{M}$.
    \item[\it{(c)}] $-1<\phi ^{'}(\overline{t})\leq 0$.
    \item[\it{(d)}] $\phi ^{2}(t)>t$ for all $t\in
    [x_{M},\overline{t}).$
\end{description}
\end{Lemma}
Proof. (a) We use the following equation in our proof frequently
\begin{equation}\label{f4}
bx_{M}^{2}+2cx_{M}+3d=0,
\end{equation}
which is clear by the fact that $x_{M}$ is an extremum of $\phi $.
Now, we show that $G_{1}$ is positive. Eq.(\ref{f4}) implies that
\begin{eqnarray}\label{f5}
  G_{1} & = & -dx_{M}^{4}+adx_{M}^{3}-cx_{M}(bx_{M}^{2}+cx_{M})-d(bx_{M}^{2}+2cx_{M}+d)
  \\
\nonumber   & = & -dx_{M}^{4}+adx_{M}^{3}+c^{2}x_{M}^{2}+3cdx_{M}+2d^{2}.
\end{eqnarray}

Now consider the polynomial $P$ in Theorem 2. Recall that the roots
of $P$ are the equilibria of $\phi$. Since $x_{M}\leq \overline{t}$
then $P(x_{M})<0$, i.e., $x_{M}^{4}<ax_{M}^{3}+bx_{M}^{2}+cx_{M}+d$.
This fact together with (\ref{f5}) imply that
\begin{equation}\label{f6}
G_{1}>(c^{2}-bd)x_{M}^{2}+2cdx_{M}+d^{2},
\end{equation}
The fact that $c<c^{*}$ together with (\ref{f4}) and (\ref{f6})
yield
\begin{equation}\label{f7}
G_{1}>2bdx_{M}^{2}+2cdx_{M}+d^{2}=d(2bx_{M}^{2}+2cx_{M}+d)=d(-2cx_{M}-5d),
\end{equation}
Since $c<c^{*}$ then the inequality $-c\sqrt{c^{2}-3bd}> -c^{2}+3bd$
holds. Some algebra show that this inequality is equivalent to the
following inequality
\begin{equation}\label{f8}
-cx_{M}> 3d,
\end{equation}
So $-2cx_{M}-5d>0$. Thus, by (\ref{f8}) $G_{1}>0$.

Next, we show that $G_{2}$ is positive. The inequality $c^{2}>3bd$,
(\ref{f4}), and (\ref{f8}) imply that
\begin{eqnarray*}
  G_{2} & = & a^{2}x_{M}^{6}+(ab-c)x_{M}^{5}-2dx_{M}^{4}+c^{2}x_{M}^{2}+ax_{M}^{3}(bx_{M}^{2}+
  2cx_{M}+2d)+
  \\
   &  &bx_{M}^{2}(bx_{M}^{2}+2cx_{M}+d)+d(bx_{M}^{2}+2cx_{M}+d) \\
   & =
   &a^{2}x_{M}^{6}+(ab-c)x_{M}^{5}-2dx_{M}^{4}-adx_{M}^{3}+(c^{2}-2bd)x_{M}^{2}-2d^{2}\\
   &
   >&a^{2}x_{M}^{6}+(-cx_{M}-2d)x_{M}^{4}+a(bx_{M}^{2}-d)x_{M}^{3}+d(bx_{M}^{2}-2d)\\
   &=&a^{2}x_{M}^{6}+(-cx_{M}-2d)x_{M}^{4}+a(-2cx_{M}-4d)x_{M}^{3}+d(-2cx_{M}-5d)\\
   &>&a^{2}x_{M}^{6}+dx_{M}^{4}+2adx_{M}^{3}+d^{2}\\
   &>&0.
\end{eqnarray*}

(b) Since $x_{M}\leq \overline{t}$ then by Theorem 3 we obtain that
$\phi (x_{M})\geq x_{M}$. Also, it's easy to verify that
$$G(x_{M})=aG_{2}+x_{M}G_{1},$$
where $G$ is as in Lemma \ref{Lemma2}. So by (a) we have
$G(x_{M})>0$. Therefore, the right hand side of (\ref{f3}) (for
$t=x_{M}$) is positive. Hence, $\phi ^{2}(x_{M})>x_{M}$.

(c) It's evident that $\phi ^{'}(\overline{t})\leq 0$. It remains to
verify that $\phi ^{'}(\overline{t})>-1$. Some calculations show
that this inequality is equivalent to the inequality $a\overline{t}^{3}-c\overline{t}-2d>0,$ which is true by (\ref{f8}) and the fact that $\overline{t}\geq
x_{M}$.

(d) Assume, for the sake of contradiction, that (d) is not true. By
Lemma 1, $\phi$ has at most three 2-cycles. It's evident that this
fact together with (b) and (c) imply that one and only one of the
following cases is possible:
\begin{description}
    \item[\it{(i)}] either $\phi $ has two prime period two
    solutions in the interval $[x_{M},\overline{t})$ or,
    \item[\it{(ii)}] $\phi $ has a unique prime period two solution
    in the interval $[x_{M},\overline{t})$ so that $\phi ^{2}$ is
    tangent to the 45-degree line at this point.
\end{description}

Assume that $p$ is the greatest prime period two solution of $\phi $
in $[x_{M},\overline{t})$. Then it's evident that $(\phi^{2})'(p)>1$ if $(i)$ holds and $(\phi ^{2})'(p)=1$ if
$(ii)$ holds. Therefore, in both cases we have
\begin{equation}\label{formula4}
(\phi ^{2})'(p)\geq 1,
\end{equation}
Consider the 2-cycle $(p,q)$. So by Eq.(\ref{formula1}) we have
\begin{equation}\label{formula5}
qp^{3}=ap^{3}+bp^{2}+cp+d, \ \ \ pq^{3}=aq^{3}+bq^{2}+cq+d.
\end{equation}

On the other hand, since $p,q>x_{M}$ then we obtain from (\ref{f8})
that $ap^{3}-cp-2d>0$ and $aq^{3}-cq-2d>0$. Therefore
\begin{equation}\label{formula6}
\frac{bp^{2}+2cp+3d}{ap^{3}+bp^{2}+cp+d}<1, \ \ \
\frac{bq^{2}+2cq+3d}{aq^{3}+bq^{2}+cq+d}<1,
\end{equation}\
By (\ref{formula5}), (\ref{formula6}), and the chain role of
calculus one can write
\begin{eqnarray*}
  (\phi ^{2})'(p)=\phi '(p)\phi
'(q) &=& \frac{bp^2+2cp+3d}{p^4}.\frac{bq^2+2cq+3d}{q^4} \\
   &=& \frac{bp^2+2cp+3d}{qp^3}.\frac{bq^2+2cq+3d}{pq^3} \\
   &=& \frac{bp^2+2cp+3d}{ap^3+bp^2+cp+d}.\frac{bq^2+2cq+3d}{aq^3+bq^2+cq+d} \\
   &<& 1,
 \end{eqnarray*}
which simply contradicts (\ref{formula4}). The proof is complete.

\begin{Lemma}\label{Lemma3} Assume that $c\in (c_{-},c^{*})$ and let $c_{1}^{*}=-2\sqrt{bd}$.
\begin{description}
    \item[\it{(a
    )}] $c_{-}<c_{1}^{*}$.
    \item[\it{(b)}] $\phi (x_{m})>a$ if and only if $c>c_{1}^{*}$.
    \item[\it{(c)}] Suppose that $c>c_{1}^{*}$. Then there exists a
    unique number $\eta >x_{M}$ such that $\phi (\eta )=\phi (x_{m})$
    where
    $$\eta =\frac{-dx_{m}}{cx_{m}+2d}.$$
\end{description}
\end{Lemma}
Proof. (a) Consider the cubic polynomial $Q$ in Theorem 1. Some
computations show that $Q(c_{1}^{*})=4abd\sqrt{bd}+27a^{2}d^{2}>0$.
This fact together with the fact that $c_{-}$ is the unique negative
root of $Q$ prove (a).

(b) The inequality $\phi (x_{m})>a$ is equivalent to the inequality
$bx_{m}^{2}+cx_{m}+d>0$ which by the equation
$bx_{m}^{2}+2cx_{m}+3d=0$ is equivalent to $-cx_{m}-2d>0$. Some
algebra show that the later inequality is equivalent to
$c\sqrt{c^{2}-3bd}>-c^{2}+2bd$. Note that both sides of this
inequality are negative since $c<c^{*}$. By squaring both sides of
this inequality and simplifying it we obtain $c^{2}<4bd$ or,
equivalently $c>c_{1}^{*}$.

(c) By (b) we have $\phi (x_{m})>a$. Also we know that $\phi
(x_{M})>\phi (x_{m})$. Therefore, by the intermediate value theorem
and the fact that $\phi $ is decreasing on $(x_{M},\infty )$ there
should exists a unique number $\eta >x_{M}$ such that $\phi (\eta
)=\phi (x_{m})$. Simplify the equation $\phi (\eta )=\phi (x_{m})$
to obtain
$$(bx_{m}^{2}+cx_{m}+d)\eta ^{2}+(cx_{m}^{2}+dx_{m})\eta +dx_{m}=0, $$
which using the equation $bx_{m}^{2}+2cx_{m}+3d=0$ is equivalent to
the following equation
$$(-cx_{m}-2d)\eta ^{2}+(cx_{m}^{2}+dx_{m})\eta +dx_{m}=0 ,$$
some algebra show that the later equation equals the following
equation
$$(x_{m}-\eta )((cx_{m}+2d)\eta +dx_{m})=0,$$
Therefore since $\eta >x_{m}$ we have $\eta
=-dx_{m}/(cx_{m}+2d)$. The proof is complete.

\begin{Lemma}\label{Lemma4} Assume that $(p,q)$ is a 2-cycle of
$\phi$ with $p<q$ and $c_{m}$ and $c_{M}$ are as in Theorem 2.
\begin{description}
\item[\it{(a)}] If $\phi $ has a unique equilibrium $\overline{t}$ then
$$p<\overline{t}<q.$$
\item[\it{(b)}] If $\phi $ has two equilibria
$\overline{t}_{1}<\overline{t}_{2}$ then for $c=c_{m}$ one of the
following cases is possible
$$p<\overline{t}_{1}<q<\overline{t}_{2},\ \ \  \text{or}\ \ \  p<\overline{t}_{1},\ q>\overline{t}_{2},$$
and for $c=c_{M}$ the following case holds
$$p<\overline{t}_{1},\ q>\overline{t}_{2}.$$
\item[\it{(c)}] If $\phi $ has three equilibria
$\overline{t}_{1}<\overline{t}_{2}<\overline{t}_{3}$ then one of the
following cases is possible
$$p<\overline{t}_{1}<q<\overline{t}_{2},\ \ \  or\ \ \  p<\overline{t}_{1},\ q>\overline{t}_{3}.$$
\end{description}
\end{Lemma}
Proof. (a) If $p<q<\overline{t}$ then by Theorem 3(a), $p=\phi
(q)>q$, a contradiction. A similar contradiction obtains when
$\overline{t}<p<q$. Lemma \ref{Lemma2}(d) plays an essential role
for the rest of theorem whose proofs are somehow easy and similar
and will be omitted.\\

The following theorem is about the convergence of solutions of
Eq.(\ref{formula1}) when $c\in [c^{*},\infty )$.
\begin{Theorem}\label{t6} Assume that $c\geq c^{*}$, $\overline{t}$ is the unique equilibrium of $\phi $,
and $\{t_{n}\}_{n=0}^{\infty } $ is a positive solution for
Eq.(\ref{formula1}). Consider polynomial $G$ in Lemma 1 and suppose
that $G$ has no iterated root of order even. Then there are four
cases to consider as follow:
\begin{description}
  \item[\it{(a)}] $\phi $ has no 2-cycle; In this case $\{t_{n}\}$ converges to $\overline{t}$.
  \item[\it{(b)}] $\phi $ has one 2-cycle $(p,q)$ with
  $p<\overline{t}<q$; In this case $\{t_{n}\}$ converges to the $2$-cycle
  $(p,q)$ if $t_{0}\neq \overline{t}$. Otherwise, $\{t_{n}\}$ simply converges to $\overline{t}$.
  \item[\it{(c)}] $\phi $ has two 2-cycles $(p_{1},q_{1})$, $(p_{2},q_{2})$ with
  $p_{1}<p_{2}<\overline{t}<q_{2}<q_{1}$; In this case $\{t_{n}\}$ converges to
  $\overline{t}$ if $t_{0}\in (p_{2},q_{2})$ and converges to the $2$-cycle $(p_{2},q_{2})$ if $t_{0}=p_{2}$ or $q_{2}$. Otherwise, $\{t_{n}\}$ converges
  to the $2$-cycle $(p_{1},q_{1})$.
  \item[\it{(d)}] $\phi $ has three 2-cycles $(p_{1},q_{1})$, $(p_{2},q_{2})$, $(p_{3},q_{3})$
  with $p_{1}<p_{2}<p_{3}<\overline{t}<q_{3}<q_{2}<q_{1}$; In this case $\{t_{n}\}$ converges to
  the $2$-cycle $(p_{3},q_{3})$ if $t_{0}\in (p_{2},q_{2})\setminus  \{\overline{t}\}$, converges to $\overline{t}$ if
  $t_{0}=\overline{t}$, and converges to the 2-cycle $(p_{2},q_{2})$ if $t_{0}=p_{2}$ or $q_{2}$. Otherwise, $\{t_{n}\}$ converges to the
  2-cycle $(p_{1},q_{1})$.
\end{description}
\end{Theorem}
Proof. $\phi ^{2}$ is increasing on $(0,\infty )$ since $\phi $ is
decreasing on $(0,\infty )$. This is a key point in this theorem and
facilitates the proof. The proof is somehow similar for all cases.
So we only give the proof for one of them. Let's prove (b). By Lemma
\ref{Lemma4}(a), $p<\overline{t}<q$. Also, $\phi ^{2}(t) \rightarrow
a$ as $t\rightarrow 0^+$ and $\phi ^{2}(t)\rightarrow \phi (a)$ as
$t\rightarrow \infty $. On the other hand, since $G$ has no iterated
root of order even then $\phi ^2$ is not tangent to the 45 degree
line at $p$ and $q$. Thus
$$(\phi ^2(t)-t)(t-p)(t-\overline{t})(t-q)<0.$$

The proof of convergence is easy using the above inequality and the
fact that $\phi ^2$ is increasing. The proof is complete.

\begin{Remark}\label{Remark1} In Theorem \ref{t6}, it is assumed (for the sake of simplicity) that
$G$ has no iterated root of order even. This assumption is not
necessary. In fact if $(p,q)$ is a 2-cycle for $\phi $ in which $p$
and $q$ are iterated roots of $G$ of order even then $\phi ^2$ is
tangent to the 45 degree line at $p$ and $q$. In this case $\phi ^2$
is semiasymptotic stable (attracting from one side and repelling
from the other side) at points $p$ and $q$. There is a similar
theorem for such a case and there's no need to mention it.

On the other hand, all cases in Theorem \ref{t6} are possible. In
all of the following examples cases (a)-(d) in Theorem 4 occur
respectively:
\begin{description}
    \item[\it{(i)}] If $a=b=c=d=1$ then
$\phi $ has no 2-cycle.
    \item[\it{(ii)}] If $a=0.1,b=2,c=1,d=0.1$ then $\phi $ has a
unique 2-cycle $C=(0.1118,169.4132)$.
    \item[\it{(iii)}] If $a=0.21,b=2.1,c=$-$2.8,d=1.3$ then $\phi $ has two 2-cycles $C_{1}=(0.2593,41.2206)$ and $C_{2}=(0.3525,13.3090)$.
    \item[\it{(iv)}] If $a=0.18,b=2.1,c=-2.8,d=1.3$ then $\phi $ has three
    2-cycles as follow\\
    $$C_{1}=(0.2001,102.9321),\ C_{2}=(0.4058,7.8071),\ C_{3}=(0.7646,1.0453).$$
\end{description}
\end{Remark}

The following theorem treats the convergence of solutions of
Eq.(\ref{formula1}) when  $c\in (c_{-},c^{*})$ and $\phi $ has a
unique equilibrium $\overline{t}$ with $x_{M}\leq
\overline{t}$ or $x_{m}\leq \overline{t}\leq x_{M}$.
\begin{Theorem}\label{Theorem5} Assume that $c\in (c_{-},c^{*})$, $\phi
$ has a unique equilibrium $\overline{t}$, $x_{M}\leq \overline{t}$
or $x_{m}\leq \overline{t}\leq x_{M}$, and $\{t_{n}\}_{n=0}^{\infty
}$ is a positive solution for Eq.(\ref{formula1}).

\begin{description}
    \item[\it{(a)}] either $\phi $ has no 2-cycle or, it has two
    2-cycles $(p_{1},q_{1}),(p_{2},q_{2})$ with $p_{1}\leq p_{2}<x_{m}<\overline{t}<q_{2}\leq
    q_{1}$ (note that if $p_{1}=p_{2}=p$ then there is essentially
    one 2-cycle and $p$ is an iterated root (of order 2) of the polynomial
    $G$ in Lemma \ref{Lemma1}).
    \item[\it{(b)}] If $\phi $ has no 2-cycle then $\{t_{n}\}$
    converges to $\overline{t}$.
    \item[\it{(c)}] Assume that $\phi $ has two 2-cycles $(p_{1},q_{1}),(p_{2},q_{2})$
    with $p_{1}\leq p_{2}<\overline{t}<q_{2}\leq q_{1}$. Then
    $\{t_{n}\}$ converges to $\overline{t}$ if $t_{0}\in
    (p_{2},q_{2})$ and converges to the 2-cycle $(p_{2},q_{2})$ if
    $t_{0}=p_{2}$ or $q_{2}$. Otherwise, $\{t_{n}\}$ converges to
    the 2-cycle $(p_{1},q_{1})$.
\end{description}
\end{Theorem}
Proof. We only give the proof for the case $x_{M}\leq \overline{t}$.
The proof of the case $x_{m}\leq \overline{t}\leq x_{M}$ is more
easier and somehow similar and will be omitted. Before proceeding to
proof note that $\phi $ is decreasing on the intervals $(0,x_{m}]$
and $[x_{M},\infty )$, and increasing on the interval
$[x_{m},x_{M}]$. Also $\phi ^2(t)\rightarrow a$ as $t\rightarrow
0^+$ and by Theorem \ref{t3}(a), $(\phi (t)-t)(t-\overline{t})<0$
for all $t>0,\ t\neq \overline{t}$. We use these facts in the proof
frequently (without mentioning them again).

(a) At first we show that $\phi ^{2}(t)>t $ on
$[x_{m},\overline{t})$. By Lemma \ref{Lemma2}(d) $\phi ^{2}(t)>t$ on
$[x_{M},\overline{t})$. So it suffices to show that $\phi ^{2}(t)>t$
on $[x_{m},x_{M})$. Assume that $t\in [x_{m},x_{M})$. Then $t< \phi
(t)<\phi (x_{M})$. So, either $t< \phi (t)\leq \overline{t}$ or,
$\overline{t}<\phi (t)<\phi (x_{M})$. In the former case Theorem
\ref{t3}(a) yields $t<\phi (t)<\phi ^{2}(t)$ and in the later case
Lemma \ref{Lemma2}(c) implies that $t<x_{M}<\phi ^{2}(x_{M})<\phi
^{2}(t)$.

Next, (by the previous discussions) if $\phi ^2(t)>t$ on $(0,x_{m})$
then $\phi $ has no 2-cycle. Otherwise, there exists $p<x_{m}$ such
that either $\phi ^2$ crosses the 45 degree line at $p$ or, it is
tangent to the 45 degree line at $p$. Let $q=\phi (p)$. In the
former case since $q$ is another period 2 point for $\phi $ and $G$
is a polynomial of degree 6 then by the intermediate value theorem
there should exists exactly another 2-cycle for $\phi $ other than
$(p,q)$. In the later case $p$ is an iterated root for $G$ and by
the same reasons we obtain that $\phi$ has no other 2-cycle.

(b) Note that in this case
\begin{equation}\label{ff1}
(\phi ^2(t)-\overline{t})(t-\overline{t})<0, \ \ \ t\neq
\overline{t}.
\end{equation}

Consider the interval $I=[x_{M},\phi (x_{M})]$. Using Lemma
\ref{Lemma2}, it's easy to verify that $I$ is invariant and $\phi
^2$ is increasing on $I$. These facts together with (\ref{ff1})
imply that $\{t_{n}\}$ converges to $\overline{t}$ if $t_{0}\in I$.
If we show that all of iterates of $\phi $ will eventually end up in
$I$ then the proof is complete. Note that it's sufficient to prove
this claim when $t_{0}\in [x_{m} ,x_{M})\cup (\phi (x_{M}),\infty
)$. Now, suppose that $t_{0}\in [x_{m},x_{M})$. Then $t_{n}\in
(x_{m},\phi (x_{M}))$ for all $n\in \Bbb{N}$. If $t_{n}\in
(x_{m},x_{M})$ for all $n\in \Bbb{N}$ then by Theorem (\ref{t3})(a),
$\{x_{n}\}$ is increasing and hence convergent to a number in
$(x_{m},x_{M})$ which simply is a contradiction. Thus $\{t_{n}\}$
eventually ends up in $I$ in this case.

On the other hand suppose that $t_{0}\in (\phi (x_{M}),\infty )$. We
claim that for some $n_{0}\in \Bbb{N}$ $t_{n_{0}}\in [x_{m},\phi
(x_{M})]$ and therefore we are done. For the sake of contradiction,
assume that such a claim is not true. Then, $t_{2n}>\phi (x_{M})$
and $t_{2n+1}<x_{m}$ for all $n\in \Bbb{N}$. So by (\ref{ff1}) it
could be shown that $\{t_{2n}\}$ is decreasing and $\{t_{2n+1}\}$ is
increasing. Hence, $\{x_{n}\}$ converges to a 2-cycle, a
contradiction.

(c) In this case we have
$$(\phi ^{2}(t)-t)(t-p_{1})(t-p_{2})(t-\overline{t})(t-q_{2})(t-q_{1})<0, \ \ \ t\neq  p_{1},p_{2},\overline{t},q_{1},q_{2}.$$

By (a) we know that $p_{1}<p_{2}<x_{m}$. Using this fact and Theorem
\ref{t3}(a) it's easy to show that $\phi ^{2}$ is increasing if
$t\in (0,p_{2})\cup (q_{2},\infty )$. Therefore, if $t_{0}\in
(0,p_{2})\cup (q_{2},\infty )$ then the proof is similar to Theorem
\ref{t6} and will be omitted. On the other hand, if $t_{0}\in
(p_{2},q_{2})$ then the proof is exactly similar to what
is applied in (b). The proof is complete.
\begin{Remark}\label{Remark2} The following examples are some
examples for all of cases in Theorem \ref{Theorem5}. Examples
$(i)$ and $(ii)$ represent cases $(b)$ and $(c)$
(respectively) when $x_{M}\leq \overline{t}$ while examples
$(iii)$ and $(iv)$ represent the same cases (respectively)
when $x_{m} \leq \overline{t}\leq x_{M}$.
\begin{description}
    \item[\it{(i)}] If $a=1,b=5,c=-4,d=1$ then $\phi $ has no
    2-cycle.
    \item[\it{(ii)}] If $a=0.1,b=5,c=-4,d=1$ then $\phi $ has two
    2-cycles $C_{1}=(0.1111,450.5876)$ and $C_{2}=(0.2019,48.2751)$.
    \item[\it{(iii)}] If $a=0.15,b=4,c=-4,d=1.1$ then $\phi $ has no
    2-cycle.
    \item[\it{(iv)}] If $a=0.1,b=4,c=-4,d=1.1$ then $\phi $ has two
    2-cycles $C_{1}=(0.1068,590.5885)$ and
    $C_{2}=(0.2378,28.0116)$.\\
\end{description}
\end{Remark}

The following theorem discusses (in some details) about the
convergence of solutions of Eq.(\ref{formula1}) when $c\in
(c_{-},c^{*})$ and $\phi $ has a unique equilibrium $\overline{t}$
with $\overline{t}<x_{m}$.
\begin{Theorem}\label{Theorem6} Assume that $c\in (c_{-},c^{*})$, $\phi $ has a unique
equilibrium $\overline{t}$ with $\overline{t}<x_{m}$, and the
sequence $\{t_{n}\}_{n=0}^{\infty }$ is a positive solution for
Eq.(\ref{formula1}). Consider the quantities $c^{*}_{1}$ and $\eta $
in Lemma 3 and let $I=[\phi (x_{m}),\phi ^{2}(x_{m})]$.
\begin{description}
    \item[\it{(a)}] Under the following hypothesis the interval $I$ is invariant.\\
      \emph{(H)} Either $c\leq c^{*}_{1}$ or, $c>c^{*}_{1}$ but
      $\phi ^{2}(x_{m})\leq \eta $.\\
Moreover, if $c\leq c^{*}_{1}$ then all iterations of $\phi $ will
eventually end up in $I$.
    \item[\it{(b)}] If $\phi $ has no 2-cycle then $\{t_{n}\}$
    converges to $\overline{t}$.
    \item[\it{(c)}] Assume that $\phi $ has one 2-cycle $(p,q)$ with $p<\overline{t}<q\leq
    x_{m}$. Let $\mathcal{S}=\{\phi
    ^{-n}(\overline{t})\}_{n=0}^{\infty }$. If $\overline{t}<-dx_{M}/(cx_{M}+2d)$
    then $\mathcal{S}=\{\overline{t}\}$. Also, if $t_{0}\not \in
    \mathcal{S}$ then $\{t_{n}\}$ converges to the 2-cycle $(p,q)$
    otherwise, it converges to $\overline{t}$.
\end{description}
\end{Theorem}
Proof. Before proceeding to proof note that Theorem (\ref{t3})(a),
monotonic properties of $\phi $, and this fact that $a$ is a
horizontal asymptote of $\phi $ are used in the proof frequently. So
we don't mention them again.

(a) Suppose that $t\in I$. We show that $\phi (t)\in I$ if (H)
holds. At first assume that $c\leq c^{*}_{1}$. So by Lemma
\ref{Lemma3}(b) $\phi (x_{m})\leq a$. Therefore, $\phi (t)>\phi
(x_{m})$ for all $t>0$. It remains to show that $\phi (t)\leq \phi
^{2}(x_{m})$. It's easy to verify that $\phi (t)\leq \max \{\phi
^2(x_{m}),x_{M}\}$. If $\max \{\phi ^2(x_{m}),x_{M}\}=\phi
^2(x_{m})$ then we are done. Thus, assume that $\phi
^2(x_{m})<x_{M}$. If $\phi ^2(x_{m})\leq x_{m}$ or $\phi
^2(x_{m})>x_{m}$ but $t<x_{m}$ then $t<x_{m}$ and hence $\phi
(t)<\phi ^2(x_{m})$. Thus, suppose that $x_{m}<t\leq \phi
^2(x_{m})$. Therefore $\phi (t)<\phi ^3(x_{m})<\phi ^2(x_{m})$.
Hence $\phi (t)\in I$.

Next assume that $c>c^{*}_{1}$ and $\phi ^{2}(x_{m})\leq \eta $.
Similar to the previous discussions it could be shown that $\phi
(t)\leq \phi ^2(x_{m})$. It remains to show that $\phi (x_{m})\leq
\phi (t)$. This matter is obvious If $t\leq x_{M}$. Thus, suppose
that $t>x_{M}$. Therefore, since $t\leq \phi ^2(x_{m})\leq \eta $
$$\phi (t)\geq \phi ^3(x_{m})\geq \phi (\eta )=\phi (x_{m}),$$
where the equality $\phi (\eta )=\phi (x_{m})$ holds by the
definition of $\eta $ in Lemma \ref{Lemma3}. This proves the
invariance of $I$.

Finally assume that $c\leq c^{*}_{1}$. By Lemma \ref{Lemma3}(b)
$\phi (x_{m})\leq a$. So $t_{n}\geq \phi (x_{m})$ for all $n\geq 1$.
As a result, if there exists $n_{0}\in \Bbb{N}$ such that
$t_{n_{0}}\leq \phi ^2(x_{m})$ then $t_{n_{0}}\in I$ and we are
done. Thus we assume that $t_{n}>\phi ^2(x_{m})$ for every $n\geq
1$. Consequently, the sequence $\{t_{n}\}$ is decreasing and hence
convergent to a number greater than or equal to $\phi ^2(x_{m})$, a
contradiction. Therefore, all iterates of $\phi $ will eventually
end up in $I$.

(b) Note that in this case (\ref{ff1}) holds and the proof is
somehow similar to Theorem \ref{Theorem5}(b) and therefore it will
be omitted.

(c) Note that if $\phi (x_{M})<\overline{t}$ then
$\mathcal{S}=\{\overline{t}\}$ obviously. Using (\ref{f4}) and some
algebra somehow similar to what is applied in Lemma 3(c) one can
write
$$\phi (x_{M})-\overline{t}=-\frac{(x_{M}-\overline{t})^{2}}{x_{M}^{3}\overline{t}^{3}}[(cx_{M}+2d)\overline{t}+dx_{M}],$$
so $\phi (x_{M})<\overline{t}$ if and only if
$(cx_{M}+2d)\overline{t}+dx_{M}>0$. By (\ref{f8}) the later
inequality is equivalent to $\overline{t}<-dx_{M}/(cx_{M}+2d)$.
Therefore, if $\overline{t}<-dx_{M}/(cx_{M}+2d)$ then
$\mathcal{S}=\{\overline{t}\}$.

On the other hand, it's easy to show that in this case the following
inequality holds
\begin{equation}\label{f10}
(\phi ^{2}(t)-t)(t-p)(t-\overline{t})(t-q)<0, \ \ \ t\neq
p,\overline{t},q.
\end{equation}

With the help of (\ref{f10}) and an analysis somehow similar to the
previous theorems the rest of proof is easy and will be omitted. The proof is complete.

\begin{Remark}\label{Remark3} In both of the following examples
hypothesis \emph{(H)} in Theorem \ref{Theorem6} holds. Also, $(i)$
and $(ii)$ represent cases (b) and (c) in Theorem \ref{Theorem6}
respectively.
\begin{description}
    \item[\it{(i)}] If $a=0.7,b=2.2,c=-3,d=1$ then $\phi $ has no
    2-cycle.
    \item[\it{(ii)}] If $a=b=1,c=-3.3,d=3$ then $\phi $ has a unique
    2-cycle $C=(1.1687,1.3190)$.\\
\end{description}
\end{Remark}

The following theorem discusses (in some details) about the
convergence of solutions of Eq.(\ref{formula1}) when $\phi $ has two
equilibria. Theorem 2 together with Lemma 2 play an essential role
for it's proof. But, it's proof will be omitted since it is somehow
similar to the proofs of some of the previous theorems of this
section. Also, a similar theorem exists when $\phi $ has three
equilibria. So it will be omitted.\\

\begin{Theorem}\label{Theorem7} Assume that $c\in (c_{-},c^{*})$, $\phi $ has two equilibria
$\overline{t}_{1}, \overline{t}_{2}$ with
$\overline{t}_{1}<\overline{t}_{2}$, and $\{t_{n}\}_{n=0}^{\infty }$
is a positive solution for Eq.(\ref{formula1}). Also consider the
values $c_{m}$ and $c_{M}$ of the parameter $c$ in Theorem \ref{t2}.
Then there are two cases as follow:
\begin{description}
    \item[\it{(a)}] $c=c_{m}$; assume that $\delta <x_{m}$ is the (unique) number such that $\phi (\delta
                    )=\overline{t}_{2}$ and $\overline{t}_{2}\leq a$.In this case the following cases are possible:
         \begin{description}
        \item[\it{($a_{1}$)}] $x_{m}\leq \overline{t}_{1}$; If $t_{0}\in (\delta
                           ,\overline{t}_{2})$ then $\{t_{n}\}$ converges to $\overline{t}_{1}$
                           otherwise, it converges to $\overline{t}_{2}$.
        \item[\it{($a_{2}$)}] $\overline{t}_{1}<x_{m}$; In this case we consider the following cases:
             \begin{description}
             \item[\it{($a_{21}$)}] $\phi $ has no 2-cycle; If $t_{0}\in (\delta
                                    ,\overline{t}_{2})$ then $\{t_{n}\}$ converges to
                                    $\overline{t}_{1}$ otherwise, it converges
                                    to $\overline{t}_{2}$.
             \item[\it{($a_{22}$)}] $\phi $ has one 2-cycle $(p,q)$
                                    with $p<\overline{t}_{1}<q\leq x_{m}$; Let $I=[\phi (x_{m}),\phi
                                    ^{2}(x_{m})]$ and $\mathcal{S}=\{\phi ^{-n}(\overline{t}_{1})\}_{n=0}^{\infty
                                   }$. Then, $\{t_{n}\}$ converges to the 2-cycle $(p,q)$
                                   if $t_{0}\in (\delta ,\overline{t}_{2})\setminus \mathcal{S}
                                   $ and converges to $\overline{t}_{1}$ if $t_{0}\in
                                   \mathcal{S}$. Otherwise, it converges to
                                   $\overline{t}_{2}$. In
                                   particular, if $t_{0}\in I\setminus \{\overline{t}_{1}\}$ then
                                   $\{t_{n}\}$ converges to the
                                   2-cycle $(p,q)$.
             \end{description}
        \end{description}
    \item[\it{(b)}] $c=c_{M}$; assume that $\delta ^{'}<x_{m}$ is
    the (unique) number such that $\phi (\delta ^{'} )=\overline{t}_{1}$ and $\overline{t}_{1}\leq
    a$. Then $\{t_{n}\}$ converges to $\overline{t}_{1}$ if $t_{0}\in [\delta
    ^{'},\overline{t}_{1}]$ otherwise, it converges to
    $\overline{t}_{2}$.
\end{description}
\end{Theorem}

\begin{Remark}\label{Remark4} In Theorem \ref{Theorem7}(a) and Theorem \ref{Theorem7}(b) it is assumed,
for the sake of simplicity, that $\overline{t}_{2}\leq a$ and
$\overline{t}_{1}\leq a$ respectively. These assumptions are not
necessary. Since similar theorems exist without these assumptions we
don't mention them. Also, the following examples represent cases
$(a_{1}),(a_{21}),(a_{22})$, and $(b)$ in Theorem 7, respectively.
\begin{description}
    \item[\it{(i)}]  $a=1,b=2.4,c=-3.8,d=1.4$.
    \item[\it{(ii)}] $a=1,b=2,c=-3,d=1$.
    \item[\it{(iii)}] $a=1,b=1.9,c=-2.8,d=.9$ with one 2-cycle
    $C=(0.5573,0.5937)$.
    \item[\it{(iv)}] $a=2,b=.5,c=-3,d=1.5$.
\end{description}
\end{Remark}

\end{document}